\documentclass{article}
\usepackage[english,russian,ukrainian]{babel}
\usepackage{latexsym,amssymb}
\usepackage[cp1251]{inputenc}
\usepackage{amsfonts,amssymb}
\usepackage{euscript}
\newcounter{theorem} 
\newcounter{lemma} 
\renewcommand{\thetheorem}{\arabic{theorem}}
\renewcommand{\thelemma}{\arabic{lemma}}
\newcommand{\theor}{\par\refstepcounter{theorem}%
{\bf Theorem \thetheorem .}\,\,}
\newcommand{\lem}{\par\refstepcounter{lemma}%
{\bf Lemma \thelemma .}\,\,}

\begin{document}

\vspace*{7mm}

\Large







\noindent {\bf  Vitalii Shpakivskyi}\\

\noindent {\bf Vieta's Formulas for Quaternionic Polynomials}

\vspace{7mm}

\textbf{Abstract.}
The paper presents analogues of some Vieta formulas for quaternionic polynomials of the form $R_n(w)=\sum\limits_{m=0}^{n}A_mx^m$.
A criterion for the sphericity of the root (zero) of a polynomial $R_n(w)$ is established.


\Large

\section{Introduction}

The results of this article were presented by the author in 2007 as a bachelor's thesis at Zhytomyr State University \cite{Shp-2007}, but have not been published until now.

The pioneer of the consideration of algebraic quaternion equations was Ivan Niven \cite{Niven-41}. In the mentioned paper, I.~Niven also proposed a method
for computing the roots of  one-side polynomials. In collaboration with S. Eilenberg, he was among the first to establish the fundamental
 theorem of algebra for quaternions \cite{Eilenberg-44}.  They proved that every quaternion polynomial $R_n(w)$ has at least one
zero. Later,  B.~Beck \cite{Beck-1979} proved that if the number of zeros of a polynomial of degree $n$ is
finite, then there are exactly $n$ zeros counting possible repetition. Beck's result was re-proven in the paper \cite{PSh} by A.~Pogorui and M.~Shapiro.
In the paper \cite{Gentili-2008}, the authors proved the
 fundamental theorem of algebra for polynomials with coefficients in the skew field of quaternions and the division algebra of octonions. 
 Polynomials over division rings were investigated by B.~Gordon and
T.~S.~Motzkin \cite{Gordon-65}.
 
In \cite{Janovska-Opfer-2010}, a new characterization of the types of zeros is provided along with an algorithm for determining all zeros, 
including their types, without relying on an iterative process requiring convergence. In \cite{Serodio-2001}, a
 method for computing the zeros of a quaternionic polynomial $R_n(w)$ is developed, which is fundamentally based on Niven's algorithm. 

In \cite{Falcao-2018}, a Weierstrass-like method was proposed for simultaneously finding all the zeros of unilateral quaternionic polynomials. 
The paper includes a convergence analysis and several numerical examples demonstrating the method's performance. 

As established in \cite{Topuridze-2001}  and \cite{PSh} (cf. \cite{Khimshiashvili-2002}), the zero set of a canonical quaternionic polynomial (that is, a polynomial with only left or only right coefficients) consists of 
  $p$ points and $s$ two-dimensional spheres, where 
  $p+2s$ \emph{does not exceed} the algebraic degree of the polynomial. The aim of N.~Topuridze's paper \cite{Topuridze-2009} was to refine and extend the results on the 
  roots of quaternionic polynomials initially presented in \cite{Topuridze-2001}. Namely, in \cite[Definition 4.1]{Topuridze-2009}
   the multiplicity of the root is correctly defined so that the "Fundamental Theorem of Algebra"\, is true. Such a definition of multiplicity is proposed, that 
    $p+2s$ is \emph{equal} to the degree of a polynomial \cite[Theorem 4.1]{Topuridze-2009}. 
    Later in this article, a similar definition of multiplicity will be proposed, but in different words.
  
The thesis \cite{Ogunmefun}  explores recent work on finding the zeros of quaternionic polynomials.  
In \cite{Dospra-2017}, using hybrid B\'{e}zout matrices, necessary and sufficient conditions were established for a quaternionic polynomial
 to possess a complex root, a spherical root, or a complex isolated root. The work  \cite{Kalantari-2013} is devoted to the study
 a method for approximating zeros of a quaternion polynomial based on the approximate
complex zeros of a real polynomial. Moreover, it is showed how to derive bounds on the norm of
zeros of a quaternion polynomial.

The paper \cite{Milovanovic-2022} focuses on examining the properties of zeros in certain special quaternionic polynomials with restricted 
coefficients, where the real and imaginary components of the coefficients satisfy specific inequalities. Additionally, the articles \cite{Mierz-Shp-2008}, \cite{Szpakowski-2009}, and 
\cite{Mierz-2008} explore various types of quadratic quaternionic equations.

Recently, two works staging the analogue of Vieta's formulas in
Clifford geometric algebras have appeared. 
In the papers \cite{Shirokov-2022,Shirokov-2024}  discussed the generalization of Vieta theorem (Vieta’s formulas) to the case of Clifford geometric algebras. 
Also it was compared the generalized Vieta’s formulas with the ordinary Vieta’s formulas for characteristic polynomial containing eigenvalues. 

In Section \ref{sect-3} of the present paper, we establish analogues of some Vieta formulas for quaternionic polynomials of the form $R_n(w)=\sum\limits_{m=0}^{n}A_mx^m$.
In Section \ref{sect-4}, a criterion for the sphericity of the root of a polynomial $R_n(w)$ is also established. 
The main research tool is the \textit{basic polynomial} introduced in paper \cite{PSh}.

\section{Preliminary}

For the real quaternion we use the notations
$$\mathbb{H}(\mathbb{R})=\{w=x_0+x_1i+x_2j+x_3k: x_0,x_1,x_2,x_3\in\mathbb{R}\},$$
where the following equalities are true for the imaginary units $i,j,k$:
$$i^2=j^2=k^2=-1,\quad ij=-ji=k,\quad jk=-kj=i,\quad ki=-ik=j.$$
The number ${\rm Sc}(w):=x_0$ is called the \textit{scalar part} of
the quaternion $w$; $x_1i+x_2j+x_3k$ is called the \textit{vector part} of $w$
 and is denoted by ${\rm Vec}(w)$; $$|w|=\sqrt{x_0^2+x_1^2+x_2^2+x_3^2}$$ is 
 called the
\textit{modulus} of $w$. Quaternion $\overline{w}:=x_0-x_1i-x_2j-x_3k$ is
 called the \emph{conjugate} with respect to quaternion $w$. \emph{Right (left)
 polynomial} \, \cite{PSh} of quaternionic variable is called polynomial of the form
\begin{equation}\label{R_n}
R_n(w):=\sum\limits_{m=0}^{n}A_mw^m \qquad \left(L_n(w):=\sum\limits_{m=0}^{n}w^mA_m
\right),
\end{equation}
where $A_m,w\in\mathbb{H(R)}$.

\textbf{Definition 1.}\cite{PSh}.
The polynomial of the complex variable
$$
F_{2n}^{*}(z):=\sum\limits_{m=0}^{2n}\Biggr(\sum\limits_{
\large
\begin{array}{c}
d+r=m,\\
0\leq d, r\leq  n\\
\end{array}
}
\Large
\overline{A}_dA_r\Biggr)z^m
$$
$$
=\overline{A}_0A_0+(\overline{A}_0A_1+\overline{A}_1A_0)z$$
\begin{equation}\label{basic-polynom}
+(\overline{A}_0A_2+\overline{A}_1A_1+\overline{A}_2A_0)z^2+\ldots +\overline{A}_nA_nz^{2n}
\end{equation}
is called a \emph{basic polynomial} for polynomial $R_n(w)$, and polynomial 
$$
G_{2n}^{*}(z):=\sum\limits_{m=0}^{2n}\Biggr(\sum\limits_{
\large
\begin{array}{c}
d+r=m,\\
0\leq d, r\leq  n\\
\end{array}
}
\Large
A_d\overline{A}_r\Biggr)z^m
$$
$$
=A_0\overline{A}_0+(A_1\overline{A}_0+A_0\overline{A}_1)z+
(A_2\overline{A}_0+A_1\overline{A}_1+A_0\overline{A}_2)z^2+\ldots +
A_n\overline{A}_nz^{2n}
$$
is called a \emph{basic polynomial} for polynomial $L_n(w)$.

The quaternion $w=x_0+{\rm Vec}(w)$ is called a 
\emph{spherical root} of the polynomial (\ref{R_n}) if any quaternion of the form $w_1:=x_0+{\rm Vec}(w_1)$ such that $|{\rm Vec}(w)|=|{\rm Vec}(w_1)|=:r$ is also a root of this polynomial.

Let us consider Corollary 5 from paper \cite{PSh}:

Given a polynomial $R_n(w)$ (or $L_n(w)$), there exist a one-to-one correspondence between pairs of the complex-conjugate roots of the basic
polynomial $F_{2n}^{*}(z)$ (or $G_{2n}^{*}(z)$) and the spheres that contain the roots of  $R_n(w)$ (or $L_n(w)$).
Moreover, the complex roots $z$ of the
 basic polynomial are related to the roots $w$ of quaternionic polynomial by the
  following equalities
\begin{equation}\label{related-roots}
 {\rm Re}\, z={\rm Sc}(w),\quad
| {\rm Im}\, z|=| {\rm Vec}(w)|.
\end{equation}
(See \cite{PSh}, Theorem 4, items (c) --- (f)).

The following theorem is true.

\theor\label{teor-1-2009} \cite{PSh}, \cite{Topuridze-2009}. \emph{The set of roots $w$ of
the quaternionic polynomial $R_n(w)$ consists of isolated points and spheres
$$\{w:{\rm Sc}(w)=x_0,\quad |{\rm Vec}(w)|=r\}.$$
Moreover, the number of the isolated points  together with the double number
of the spheres is less or equal to $n$.}

\section{Multiplicity of the root}

Next, we will define multiplicity. Since the concept of multiplicity is clearly defined for complex polynomials,
 we will use the basic polynomial $F_{2n}^{*}(z)$  to define the multiplicity of the root of a quaternion polynomial $R_n(w)$.
A similar approach to the definition of multiplicity is taken in paper \cite[Definition 4.1]{Topuridze-2009}.

\textbf{Definition 2.}
  \textit{The multiplicity of a root (isolated or spherical)} $w$ of the polynomial
$R_n(w)$ is called the multiplicity of the corresponding
  complex root $z$ (or $\bar{z}$) of the basic polynomial $F_{2n}^{*}(z)$. 
  It is clear that if $z$ is a real number, then $z=\bar{z}$ and therefore when we calculate the multiplicity, we need to divide the multiplicity of the real root $z$ by two.
 \vskip2mm
  
\textbf{Remark 1.} The basic polynomial $F_{2n}^{*}(z)$ has exactly $2n$ roots, taking into account their multiplicity.
 Since the coefficients of the basic polynomial are real numbers, then in addition to the roots $z_1$, $z_2$,...,$z_n$, the roots also include the numbers 
 $\overline{z}_1$, $\overline{z}_2$,...,$\overline{z}_n$. Therefore, in Definition 2, it is necessary to calculate the multiplicity of either $z$ or $\overline{z}$.
  
    \vskip2mm
  
 Thus, the polynomial $R_n(w)$ (and $L_n(w)$) has exactly $n$ roots, taking into account multiplicity.
  Compare with Theorem 4.1 in \cite{Topuridze-2009}.
   Therefore, the following definition is correct.

\textbf{Definition 3.}
\emph{The set of roots} of polynomial $R_n(w)$ is called the set  $\{w_1,w_2,\ldots,w_n\}$  which contains each root as many times as its multiplicity. 
 If the spherical root $w$ has multiplicity $k$, then in \emph{the set of roots} we include arbitrary $k$ quaternions from this sphere.

 Let us consider a few examples.
 \vskip2mm
 
\textbf{Example 1.} \cite[Example 2.5]{Huang-2002}. Consider the quadratic quaternionic equation $R_2(w)=w^2+1 = 0$. This equation has unique spherical root
 $w_0=x_1i+x_2j+x_3k$, where $x_1^2+x_2^2+x_3^2=1$. For this equation, the basic polynomial has the form $F_{4}^{*}(z)=(z^2+1)^2$. Polynomial $F_{2n}^{*}$
  has two roots of multiplicity $2$: $z=i$ and $\overline{z}=-i$. Therefore, since the multiplicity of the root $z$ is 2, the quaternionic spherical root $w_0$
   also has a multiplicity of 2, which coincides with the degree of the polynomial. In addition, let us choose, for example, the following \emph{set of roots} 
   $\{i;j\}$.
 
  \vskip2mm
 \textbf{Example 2.} \cite[Example 2.7]{Huang-2002}. Consider the quadratic quaternionic equation $R_2(w)=(w-1)^2 = 0$. This equation has the unique solution is $w_0=1$.
  For this equation, the basic polynomial has the form $F_{4}^{*}(z)=(z-1)^4$. Polynomial $F_{2n}^{*}$
  has one \emph{real} roots of multiplicity $4$: $z=\overline{z}=1$. Therefore, since the multiplicity of the \textit{real} root $z$ is $4$, the quaternionic root $w_0$
    has a multiplicity of $4/2=2$, which coincides with the degree of the polynomial.
  So we have the following \textit{set of roots}
   $\{1;1\}$.
 
 \vskip2mm
 \textbf{Example 3.} \cite[Example 2.8]{Huang-2002}. Consider the quadratic quaternionic equation $R_2(w)=w^2+iw+\frac{1}{2}\,j = 0$. 
 This equation has the unique isolated solution is $w_0=\frac{k-i}{2}$.
  For this equation, the basic polynomial has the form $F_{4}^{*}(z)=(z^2+\frac{1}{2})^2$. Polynomial $F_{2n}^{*}$
  has two roots of multiplicity $2$: $z=\frac{i}{\sqrt{2}}$ and $\overline{z}=-\frac{i}{\sqrt{2}}$. 
  So, since the multiplicity of the root $z$ is 2, the quaternionic isolated root $w_0$
   also has a multiplicity of 2, which coincides with the degree of the polynomial. For this equation we have the following \emph{set of roots} 
   $\left\{\frac{k-i}{2};\frac{k-i}{2}\right\}$.
 
  \vskip2mm
 \textbf{Example 4.} \cite[p. 389]{PSh}. Consider the quaternionic equation $R_4(w)=(w^2+1)^2 = 0$. 
  This equation has unique spherical root
 $w_0=x_1i+x_2j+x_3k$, where $x_1^2+x_2^2+x_3^2=1$. For this equation, the basic polynomial has the form $F_{8}^{*}(z)=(z^2+1)^4$. Polynomial $F_{2n}^{*}$
  has two roots of multiplicity $4$: $z=i$ and $\overline{z}=-i$. So, since the multiplicity of the root $z$ is 4, the quaternionic spherical root $w_0$
   also has a multiplicity of 4, which coincides with the degree of the polynomial. In addition, let us choose, for example, the following \emph{set of roots} 
   $\left\{i,j,k,\frac{1}{2}\,i+\frac{\sqrt{3}}{2}\,j\right\}$.
   
   \vskip2mm
 \textbf{Example 5.} \cite[p. 389]{PSh}. Consider the quaternionic equation $$R_5(w)=(w^3+iw^2+jw+k)(w^2+1) = 0.$$ 
  This equation has the following roots $w_1=\frac{1}{\sqrt{2}}(1-j)$,  $w_2=-\frac{1}{\sqrt{2}}(1-j)$, $w_3=\{x_1i+x_2j+x_3k:x_1^2+x_2^2+x_3^2=1\}$.
  For this equation, the basic polynomial has the form $$F_{10}^{*}(z)=(z^6+z^4+z^2+1)(z^2+1)^2.$$ Polynomial $F_{10}^{*}$
  has the following single roots $z_{1,2}=e^{\pm\frac{\pi i}{4}}$, $z_{3,4}=e^{\pm\frac{5\pi i}{4}}$, and   
   two roots of multiplicity $3$: $z_{5,6,7}=i$ and $z_{8,9,10}=-i$. Therefore,  the roots $w_1$ and $w_2$ are single roots, and the sphere  $w_3=\{x_1i+x_2j+x_3k:x_1^2+x_2^2+x_3^2=1\}$
   has a multiplicity of 3. Thus, the number of roots, taking into account the multiplicity, coincides with the degree of the polynomial.
   For this equation we can chouse, for example, the following \emph{set of roots} 
   $\{w_1,w_2,i,j,k\}$.
   
   \vskip2mm
 Taking into account Definitions 2 and 3, we can refine Theorem \ref{teor-1-2009} (\cite{PSh}, \cite{Topuridze-2009}). 
 
 \vskip2mm
 \theor\label{teor-1+1-2009}
 \emph{The set of roots of
the quaternionic polynomial $R_n(w)$ (or $L_n(w)$) consists of isolated points and spheres. 
Moreover, the number of the isolated points and the spheres, taking into account their multiplicity, is equal to $n$.}

 \section{Analogues of some Viete formulas} \label{sect-3}

 Now let us consider some of Vieta's formulas.
 
The \textit{scalar product} of quaternions is defined as follows $(w,\vartheta):=\sum\limits_{m=0}^{3}x_my_m$, where $w=x_0+x_1i+x_2j+x_3k$,\,\,$\vartheta=y_0+y_1i+y_2j+y_3k$.
\vskip2mm

\lem \label{lem-1-2009} \emph{For any $w,\vartheta\in\mathbb{H(R)}$ the following
 equality is true}
$$w\overline{\vartheta}+\vartheta\overline{w}=\overline{w}\vartheta+\overline{\vartheta}w
=2(w,\vartheta).$$

\textbf{Proof.} We first convert the expression
$$w\overline{\vartheta}+\vartheta\overline{w}=\left(x_0+{\rm Vec}(w)\right)\left(y_0-
{\rm Vec}(\vartheta)\right)+\left(y_0+{\rm Vec}(\vartheta)\right)\left(x_0-
{\rm Vec}(w)\right)$$
$$=2x_0y_0-{\rm Vec}(w){\rm Vec}(\vartheta)-{\rm Vec}(\vartheta){\rm Vec}(w)=
\overline{w}\vartheta+\overline{\vartheta}w.$$
Then, we consider the sum
$${\rm Vec}(w){\rm Vec}(\vartheta)+{\rm Vec}(\vartheta){\rm Vec}(w)=
(-x_1y_1+x_1y_2k-x_1y_3j-x_2y_1k-x_2y_2$$
$$+x_2y_3i+x_3y_1j-x_3y_2i-x_3y_3)+(-y_1x_1+y_1a_2k-y_1x_3j-y_2x_1k-x_2y_2$$
$$+y_2x_3i+y_3x_1j-y_3x_2i-y_3x_3)=-2x_1y_1-2x_2y_2-2x_3y_3.$$
Using the two previous equalities, we obtain the statement of lemma.
\vskip2mm

\theor\label{teor-2-2009}\emph{For any polynomial $R_n(w)$ and for the 
set of its roots, taking into account their multiplicity, the following equality is true}
\begin{equation}\label{teor-2-2009-1}
\prod\limits_{m=1}^{n}|w_m|=\frac{|A_0|}{|A_n|}.
\end{equation}

\textbf{Proof.} Let $z_1,z_2,\ldots,z_{2n}$ be the roots of basic polynomial.
By definition (\ref{basic-polynom})
$$F_{2n}^{*}(z):=|A_n|^2z^{2n}+\ldots+(\overline{A}_1A_0+\overline{A}_0A_1)z+|A_0|^2.$$
Since the degree of the basic polynomial is $2n$, then it has exactly $2n$ roots. Let
$$z_2=\overline{z}_1, \qquad z_4=\overline{z}_3, \qquad z_{2n}=\overline{z}_{2n-1}.$$
Using the Vieta formula for the polynomial
 $F_{2n}^{*}$, we obtain
$$\prod\limits_{m=1}^{2n}z_m=\frac{|A_0|^2}{|A_n|^2}.$$
Let
$$z_1=\alpha_1+i\beta_1,\,z_2=\overline{z}_1=\alpha_1-i\beta_1,\,z_3=\alpha_3+i\beta_3,\,z_4=\overline{z}_3=\alpha_3-i\beta_3,\,\ldots,\,$$
$$z_{2n-1}=\alpha_{2n-1}+i\beta_{2n-1},\,z_{2n}=\overline{z}_{2n-1}=\alpha_{2n-1}-i\beta_{2n-1}.$$
Then
\begin{equation}\label{alpha+ibeta}
(\alpha_1+i\beta_1)(\alpha_1-i\beta_1)\cdot\ldots\cdot(\alpha_{2n-1}+i\beta_{2n-1})(\alpha_{2n-1}-i\beta_{2n-1})=\frac{|A_0|^2}{|A_n|^2},
\end{equation}
whence
\begin{equation}\label{teor-2-2009-2}
(\alpha_1^2+\beta_1^2)(\alpha_3^2+\beta_3^2)\cdot\ldots\cdot(\alpha_{2n-1}^2+\beta_{2n-1}^2)=\frac{|A_0|^2}{|A_n|^2}.
\end{equation}
Using  system (\ref{related-roots}), we obtain the equalities
$$ 
\alpha^2={\rm Sc}^2(w),\quad
|\beta|^2=|{\rm Vec}(w)|^2.
$$
Then  equality (\ref{teor-2-2009-2}) can be rewritten in the form
$$\left({\rm Sc}^2(w_1)+|{\rm Vec}(w_1)|^2\right)\left({\rm Sc}^2(w_2)+
|{\rm Vec}(w_2)|^2\right) $$
$$\ldots\times\left({\rm Sc}^2(w_n)+
|{\rm Vec}(w_n)|^2\right)=\frac{|A_0|^2}{|A_n|^2}$$
and
$$|w_1|^2\cdot|w_2|^2\cdot\ldots\cdot|w_n|^2=\frac{|A_0|^2}{|A_n|^2}.$$
The theorem is proved.

\textbf{Corollary 1.} \emph{If the polynomial $R_n(w)$ is monic (i.e., if $A_n=1$), then 
$$\prod\limits_{m=1}^{n}|w_m|=|A_0|.$$
}

\theor\label{teor-3-2009}  \emph{Let the conditions of Theorem} \ref{teor-2-2009} \emph{ are satisfied. Then }
$$\sum\limits_{m=1}^{n}{\rm Sc}(w_m)=-\frac{(A_n,A_{n-1})}{|A_n|^2}.$$

\textbf{Proof.} Let us write the basic polynomial in the form
$$F_{2n}^{*}(z)=\sum\limits_{m=0}^{2n}B_mz^m.$$
Let $z_2=\overline{z}_1,\,z_4=\overline{z}_3,\,\ldots,\,z_{2n}=\overline{z}_{2n-1}$ be 
the roots of basic polynomial. Since $z+\overline{z}=2{\rm Re}\,z$,
 then using Viete formula  for the 
 polynomial $F_{2n}^{*}$, we obtain
$$\sum\limits_{m=1}^{n}z_m=-\frac{B_{2n-1}}{B_{2n}}$$
hence
\begin{equation}\label{teor-3-2009=}
2({\rm Re}\,z_1+{\rm Re}\,z_3+\ldots+{\rm Re}\,z_{2n-1})=-\frac{B_{2n-1}}{B_{2n}}.
\end{equation}

Since for two complex conjugate roots of the basic polynomial 
corresponds one quaternionic root and equality (\ref{related-roots}) is true,
  equality (\ref{teor-3-2009=}) can be  rewritten in the form
$$2\sum\limits_{m=1}^{n}{\rm Sc}(w_m)=-\frac{B_{2n-1}}{B_{2n}},$$
where $B_{2n-1}=\overline{A}_{n-1}A_n+\overline{A}_nA_{n-1},\,B_{2n}=|A_n|^2.$ 
By the Lemma \ref{lem-1-2009} 
$$B_{2n-1}=\overline{A}_{n-1}A_n+\overline{A}_nA_{n-1}=2(A_n,A_{n-1}).$$
Then
$$2\sum\limits_{m=1}^{n}{\rm Sc}(w_m)=-\frac{2(A_n,A_{n-1})}{|A_n|^2}.$$
The theorem is proved.

\textbf{Corollary 2.} \emph{If the polynomial $R_n(w)$ is monic (i.e., if $A_n=1$), then 
$$\sum\limits_{m=1}^{n}{\rm Sc}(w_m)=-(A_n,A_{n-1}).$$
}

\lem\label{lem-2-2009}\emph{Let the conditions of Theorem} \ref{teor-2-2009} \emph{ are satisfied}
 \emph{and $A_0\neq 0.$ Then none of the roots
 of the corresponding basic polynomial is not zero.}

\textbf{Proof.} Suppose the contrary, i.~e. let exists a number $z=0$ such that
$$F_{2n}^{*}(0)=|A_n|^2\cdot0^{2n}+\ldots+(\overline{A}_1A_0+\overline{A}_0A_1)\cdot0+|A_0|^2=0.$$
Since $A_0\neq 0$, then $|A_0|^2\neq 0$. This contradiction proves the Lemma.

\theor\label{teor-4-2009}\emph{Let the conditions of Theorem} \ref{teor-2-2009} \emph{ are satisfied} \emph{and $A_0\neq 0.$ Then the following equation is true}
$$\sum\limits_{m=1}^{n}\frac{{\rm Sc}(w_m)}{|w_m|^2}=-\frac{(A_1,A_0)}{|A_0|^2}.$$

\textbf{Proof.} Consider the Vieta equality
$$z_1z_2\cdot\ldots\cdot z_{2n-1}+z_1z_2\cdot\ldots\cdot z_{2n-2}z_{2n}+\ldots +z_1z_3\cdot\ldots\cdot z_{2n}+z_2z_3\cdot\ldots\cdot z_{2n}=-\frac{B_1}{B_{2n}},$$
where $B_1,\,B_{2n}$ are corresponding coefficients of the basic polynomial. 
Using equality (\ref{alpha+ibeta}) and relation $z\cdot\overline{z}=|z|^2$, we obtain
$$|z_1|^2\cdot|z_3|^2\cdot\ldots\cdot|z_{2n-3}|^2(\alpha_{2n-1}+i\beta_{2n-1})+
|z_1|^2\cdot|z_3|^2\cdot\ldots\cdot|z_{2n-3}|^2(\alpha_{2n-1}-i\beta_{2n-1})$$
$$\ldots+|z_3|^2\cdot|z_5|^2\cdot\ldots\cdot|z_{2n-1}|^2(\alpha_1+i\beta_1)+
|z_3|^2\cdot|z_5|^2\cdot\ldots\cdot|z_{2n-1}|^2(\alpha_1-i\beta_1)=
-\frac{B_1}{B_{2n}}$$
or in equivalent form 
$$2\alpha_{2n-1}|z_1|^2\cdot|z_3|^2\cdot\ldots\cdot|z_{2n-3}|^2+\ldots+2\alpha_1|z_3|^2\cdot|z_5|^2\cdot\ldots\cdot|z_{2n-1}|^2=-\frac{B_1}{B_{2n}}.$$
Then, let us multiply and divide each term of the left-hight side on 
 the corresponding squares of modules:
$$\frac{2\alpha_{2n-1}\prod\limits_{m=1}^{n}|z_{2m-1}|^2}{|z_{2n-1}|^2}+\ldots+\frac{2\alpha_1\prod\limits_{m=1}^{n}|z_{2m-1}|^2}{|z_1|^2}=-\frac{B_1}{B_{2n}}$$
or
$$2\prod\limits_{m=1}^{n}|z_{2m-1}|^2\left(\sum\limits_{m=1}^{n}\frac{\alpha_{2m-1}}{|z_{2m-1}|^2}\right)=-\frac{B_1}{B_{2n}}.$$
Taking into account the equalities $B_{2n}=|A_{n}|^2$,\,$B_1=\overline{A}_1A_0+\overline{A}_0A_1=2(A_1,A_0)$,
 and relations (\ref{related-roots}), we obtain
$$2\prod\limits_{m=1}^{n}|w_m|^2\left(\sum\limits_{m=1}^{n}
\frac{{\rm Sc}(w_m)}{|w_m|^2}\right)=-\frac{2(A_1,A_0)}{|A_n|^2}.$$
Using Theorem \ref{teor-2-2009}, we have
$$2\frac{|A_0|^2}{|A_n|^2}\left(\sum\limits_{m=1}^{n}\frac{{\rm Sc}(w_m)}{|w_m|
^2}\right)=-\frac{2(A_1,A_0)}{|A_n|^2}.$$
Multiplying the previous equality on the $|A_n|^2$, we prove the theorem.

\textbf{Corollary 3.} \emph{If all the roots of a quaternion polynomial $R_n(x)$
 are pure vectors, then $(A_{n},A_{n-1})=0.$}

\textbf{Proof.} By the condition of corollary we have
 $\sum\limits_{m=1}^{n}{\rm Sc}(w_m)=0.$ Therefore
$$-\frac{(A_n,A_{n-1})}{|A_n|^2}=0.$$
Since $|A_n|^2\neq 0$, then $(A_n,A_{n-1})=0.$

The following corollary is result of Theorem \ref{teor-4-2009}.

\textbf{Corollary 4.} \emph{If all the roots of a quaternion polynomial $R_n(x)$
 are pure vectors, then $(A_{1},A_{0})=0.$}

\textbf{Remark 2.}  It is easy to verify that Theorems \ref{teor-2-2009}, \ref{teor-3-2009}, and \ref{teor-4-2009} are valid for the equations considered in Examples 1 --- 5.

\section{Structure properties of roots of the polynomials} \label{sect-4}

For a more detailed investigation of properties of the roots of quaternionic 
polynomials, we 
will find a general formula to calculate the degree of quaternion.

Consider the equalities:
$$w^0=1,$$
$$w^1=w,$$
$$w^2=2x_0w-|w|^2,$$
$$w^3=w^2w=(2x_0w-|w|^2)w=2x_0w^2-|w|^2w=2x_0(2x_0w-|w|^2)-|w|^2w=$$
$$=4x_0^2w-2x_0^2|w|^2-|w|^2w=(4x_0^2-|w|^2)w-2x_0|w|^2.$$
Similarly
$$w^4=(8x_0^3-4x_0|w|^2)w-(4x_0^2-|w|^2)|w|^2,$$
$$w^5=(16x_0^4-12x_0^2|w|^2-|w|^4)w-(8x_0^3-4x_0|w|^2)|w|^2.$$
The coefficients near $w$ and $|w|^2$ are real numbers, which we denote respectively by $Q_n(x_0,|w|^2)$ and $P_n(x_0,|w|^2)$. Then we get the equality
$$w^n=Q_n(x_0,|w|^2)w-P_n(x_0,|w|^2)|w|^2.$$
The functions $Q_n$ and $P_n$ are related by the recurrence formulae.

\textbf{Property 1.} \cite{PSh}. \emph{The following equalities are true}
$$
\begin{array}{l}
\displaystyle 1)\,\,\,Q_n(w)=2x_0Q_{n-1}(w)-P_{n-1}(w)|w|^2\,,\\[5mm]
\displaystyle 2)\,\,\, P_n(w)=Q_{n-1}(w)\,,\\[5mm]
\end{array}$$
\emph{where} $Q_0(w):=0$,\,$P_0(w):=-\frac{1}{|w|^2}$,\,$Q_1(w):=1$,\,$P_1(w):=0$.

\vskip2mm
The next theorem is a \textit{criterion of sphericity of root.}
\vskip2mm

\theor\label{teor-5-2009}  \emph{The quaternion $w_0$ is a spherical root of polynomial
 $R_n(w)$ if and only if}
\begin{equation}\label{teor-5-2009-0}
\sum\limits_{m=1}^{n}A_mQ_m(w_0)=0.
\end{equation}

\textbf{Proof.} \emph{Necessity.} If $w_0$ is a root of $R_n$, then $R_n(w_0)=0$, i.~e. $\sum\limits_{m=0}^{n}A_mw_0^m=0.$ Then
$$A_n\left(Q_n(w_0)w_0-P_n(w_0)|w_0|^2\right)+A_{n-1}
\left(Q_{n-1}(w_0)w_0-P_{n-1}(w_0)|w_0|^2\right)$$
$$\ldots+A_1\left(Q_1(w_0)w_0-P_1(w_0)|w_0|^2\right)+A_0=0,$$
which implies
\begin{equation}\label{teor-5-2009-1}
\left(\sum\limits_{m=1}^{n}A_mQ_m(w_0)\right)w_0-
|w_0|^2\sum\limits_{m=1}^{n}A_mP_m(w_0)+A_0=0.
\end{equation}
Since $w_0$ is a spherical root, then any quaternion
$$w_1=x_0+{\rm Vec}(w_1) \quad{\rm with}\quad 
|{\rm Vec}(w_0)|=|{\rm Vec}(w_1)|$$
is also a root of this polynomial. Therefore
$$A_n\left(Q_n(w_1)w_0-P_n(w_1)|w_1|^2\right)+A_{n-1}\left(Q_{n-1}(w_1)w_1
-P_{n-1}(w_1)|w_1|^2\right)$$
\begin{equation}\label{teor-5-2009-2}
\ldots+A_1\left(Q_1(w_1)w_1-P_1(w_1)|w_1|^2\right)+A_0=0.
\end{equation}
Denoted by
$$Q_n(w_0):=Q_n\left(x_0,|w_0|^2\right), \qquad Q_n(w_1):=Q_n\left(x_0,|w_1|^2\right),$$
$$P_n(w_0):=P_n\left(x_0,|w_0|^2\right), \qquad P_n(w_1):=P_n\left(x_0,|w_1|^2\right).$$
Since $|w_0|^2=|w_1|^2$, we obtain
\begin{equation}\label{2009----1}
Q_n(w_0)=Q_n(w_1), \qquad P_n(w_0)=P_n(w_1).
\end{equation}
Then equality (\ref{teor-5-2009-2}) takes the form
$$A_n\left(Q_n(w_0)w_1-P_n(w_0)|w_0|^2\right)+A_{n-1}
\left(Q_{n-1}(w_0)w_1-P_{n-1}(w_0)|w_0|^2\right)$$
$$\ldots+A_1\left(Q_1(w_0)w_1-P_1(w_0)|w_0|^2\right)+A_0=0$$
or equivalent
\begin{equation}\label{teor-5-2009-4}
\left(\sum\limits_{m=1}^{n}A_mQ_m(w_0)\right)w_1
-|w_0|^2\sum\limits_{m=1}^{n}A_mP_m(w_0)+A_0=0.
\end{equation}
The equality (\ref{teor-5-2009-4}) subtract from equality (\ref{teor-5-2009-1}). 
Then
$$\left(\sum\limits_{m=1}^{n}A_mQ_m(w_0)\right)(w_0-w_1)=0,$$
whence
$$\sum\limits_{m=1}^{n}A_mQ_m(w_0)=0.$$

\emph{Sufficiency.} 
Now we prove that any quaternion
$$w_1=x_0+{\rm Vec}(w_1) \quad {\rm such \,\, that}\quad|{\rm Vec}(w_0)|=
|{\rm Vec}(w_1)|$$
is a root of polynomial $R_n(w).$ A consequence of equality 
(\ref{teor-5-2009-0}) 
 is identity
\begin{equation}\label{teor-5-1-2009}
\left(\sum\limits_{m=1}^{n}A_mQ_m(w_0)\right)w_0=0.
\end{equation}
The relations (\ref{2009----1}), (\ref{teor-5-1-2009}) implies the equality
\begin{equation}\label{teor-5-2009-5}
\left(\sum\limits_{m=1}^{n}A_mQ_m(w_1)\right)w_1=0.
\end{equation}
Now, from (\ref{teor-5-1-2009}) and (\ref{teor-5-2009-5}), we obtain
\begin{equation}\label{teor-5-2009-5-0}
\left(\sum\limits_{m=1}^{n}A_mQ_m(w_0)\right)w_0=
\left(\sum\limits_{m=1}^{n}A_mQ_m(w_1)\right)w_1.
\end{equation}

Using equality (\ref{teor-5-2009-5-0}),  formula (\ref{teor-5-2009-1}) 
takes the form
$$\left(\sum\limits_{m=1}^{n}A_mQ_m(w_0)\right)w_0-
|w_0|^2\sum\limits_{m=1}^{n}A_mP_m(w_0)+A_0=0,$$
which is equivalent to the following condition
$$A_n\left(Q_n(w_1)w_1-P_n(w_1)|w_1|^2\right)+A_{n-1}\left(Q_{n-1}(w_1)w_1
-P_{n-1}(w_1)|w_1|^2\right)$$
$$\ldots+A_1\left(Q_1(w_1)w_1-P_1(w_1)|w_1|^2\right)+A_0=0,$$
that is
$$\sum\limits_{m=0}^{n}A_mw_1^m=0.$$
Therefore, $w_1$ is a root of the polynomial $R_n(w)$. The theorem is proved.

\textbf{Remark 3.}  We note that another criterion for the sphericity of the root is established in paper \cite[Proposition 2]{PSh}.

The following statement follows from Theorem \ref{teor-5-2009}.

\textbf{Corollary 5.} \emph{If $w_0$ and $\overline{w}_0$ are the 
roots of the polynomial $R_n(w)$ then this root is spherical.}

\textbf{Remark 4.} Corollary 5 is also established in paper \cite[ Corollary 2.1]{PSh}.

The following theorem establishes a relation between the roots of 
polynomials $R_n(w)$ and $L_n(w)$.

\theor\label{teor-6-2009} \emph{The following equality is true}
$$F^*_{2n}(w)=G^*_{2n}(w).$$

\textbf{Proof.} It is obviously that $\deg F^*_{2n}(w)=\deg G^*_{2n}(w)$. 
Now we prove the equality of the corresponding coefficients. For example, 
consider the second coefficients of basic polynomials corresponding 
for the right and left polynomials:
$$B_{1}^r=\overline{A}_0A_1+\overline{A}_1A_0,\qquad B_{1}^{\ell}=A_1\overline{A}_0+A_0\overline{A}_1.$$
By Lemma \ref{lem-1-2009} \, $B_{1}^r=B_{1}^{\ell}$.
 For the same reason all other coefficients also equal. The theorem is proved.

\textbf{Remark 5.}
By Theorem \ref{teor-6-2009}, the all results for the right polynomials $R_n(w)$
 also true for the
left polynomials $L_n(w)$.

 \subsection*{Acknowledgment}
1. This work was supported by a grant from the Simons Foundation
(1030291,V.S.Sh.).

2. The author expresses his sincere gratitude to the reviewers, whose comments significantly improved the paper.

\renewcommand{\refname}{References}

\vskip10mm
\noindent Vitalii Shpakivskyi

\noindent Institute of Mathematics of the\\
National Academy of Sciences of Ukraine, Kyiv

\noindent e-mail: shpakivskyi86@gmail.com

\end{document}